\documentclass{article}

\usepackage{arxiv}

\usepackage[utf8]{inputenc} 
\usepackage[T1]{fontenc}    
\usepackage{hyperref}       
\usepackage{url}            
\usepackage{booktabs}       
\usepackage{amsfonts}       
\usepackage{nicefrac}       
\usepackage{microtype}      
\usepackage{lipsum}		
\usepackage{graphicx}
\usepackage{natbib}
\usepackage{doi}

\usepackage{algorithm}%
\usepackage{algorithmicx}%
\usepackage{algpseudocode}%
\usepackage[figuresright]{rotating}%
\usepackage{graphicx}
\usepackage{mathtools}
\usepackage{url}
\usepackage{anyfontsize}
\usepackage{lmodern}

\title{An Enhanced Approach for the Dial a ride problem with drivers preferences}

\author{ Sana Ouasaid\\
	Department of Computer Science\\
	Hassan II University\\
	LIM@II-FSTM, B.P. 146, Mohammedia, Morocco\\
	\texttt{ouasaidsanaa@gmail.com} \\
	\And
	Mohammed Saddoune\\
	Department of Computer Science\\
	Hassan II University\\
	LIM@II-FSTM, B.P. 146, Mohammedia, Morocco\\
	}

\date{}

\renewcommand{\shorttitle}{E-ILS for DARPDP}

\hypersetup{
pdftitle={An Enhanced Approach for the Dial A Ride Problem with Drivers Preferences},
pdfsubject={math.OC, cs.DS},
pdfauthor={Sana Ouasaid, Mohammed Saddoune},
pdfkeywords={Dial-A-Ride Problem, Driver preferences, Metaheuristic Approach, Multi-objective optimization},
}

\begin{document}
\maketitle

\begin{abstract}
	The paper addresses a variant of the Dial-A-Ride problem with additional features. It is referred to as the DARP with driver preferences, which attempts to determine a solution more driver-oriented by designing a short trip in a specific direction that has to be finished at a destination of interest within a restricted time window.  For this purpose, two solutions are considered. The first involves solving the new MILP exactly using the CPLEX software. The second is a new approach that employs an iterated local search as a general framework and exploits many heuristics. Numerical experiments indicate that the approach can efficiently solve the generated DARPDP instances in a reasonable time.
\end{abstract}


\section*{Introduction}
\quad With the ever-growing demand for mobility, several shared services have merged with the purpose of solving specific problems. As an example, the Dial-A-Ride problem (DARP) is considered an NP-Hard problem concerning on-demand transportation. The DARP consists of planning vehicle routes to efficiently serve transport requests in different locations. Each request specifies the number of passengers and the locations for the start and the end of the trip. As well as the time windows associated with those locations. This problem arises in several contexts, such as commercial taxi services, which are defined through some key choices: fleet, passengers, route, schedule, and fares. In this paper, we propose an extension of DARP from the point of view of taxi drivers.

\quad The purpose behind this is that the drivers prefer short- and medium-distance trips on more profitable lines. The lines can be considered as connecting two target designated areas known for high demand in a specific period, such as guys, streetcar stations, etc. Indeed, starting from his current position, the taxi driver aims to pick up and drop off as many as possible requests in a flexible manner and ends his trip at a destination of interest at a specific time. Our proposal is to assist drivers by providing a solution that can be truly designed around their preferences. It will adapt routes and schedules to pick up passengers in a more profitable way for drivers while maintaining the preferences of the onboard passengers.
	
\quad The remainder of this work can be summarized as follows. Section \ref{Sec2} presents some literature reviews on the DARP variants and solution methods. Afterward, a description with related assumptions to the problem studied is given in Section \ref{Sec3}. In Section \ref{Sec4}, the general solution approach is explained in detail. Section \ref{Sec5} describes how new instances based on real-world data are generated and presents the computational results and analysis of the performance of the CPLEX solver and the proposed method. Section \ref{Sec6} concludes the paper and presents some possible future directions.

\section{Literature Review}\label{Sec2}
	
\quad Depending on the real-world applications, a DARP was extended to different variants defined with different constraints and objectives. To solve these variants several solution methods are designed. In the following, we emphasize more recent contributions, significant problem extensions, and methods applied.
 
\quad There are many variants of DARP, obtained either by adding constraints. Like the DARP with time windows (DARPTW for Dial a ride problem with time windows) in which the visit of each customer must be completed in a given time interval. Or through modification, for example, the Multi-depot DARP where a depot is specified for each vehicle rather than a single depot serving as the departure and destination of all the vehicles. For a detailed presentation of the different variants of DARP, as well as the existing solution methods. One can refer to Cordeau and Laporte \cite{bib1}.

\quad Parragh et al.\cite{bib2} tackled the DARP with split requests and profits (DARPSRP). They consider that each request concerns several persons and can be split and assigned to different vehicles. In fact, multiple requests can be handled at the same time. The problem concept is the design of routes aiming for the maximization of the total profit under a set of constraints vehicle capacity, customer time windows, maximum ride time, precedence, and pairing. 

\quad Masmoudi et al.\cite{bib3} presented a Dial-a-Ride Problem with Electric Vehicles and battery swapping stations (DARP-EV). They consider multiple types of EVs such that each EV can provide several types of resources (each with a specific capacity). The problem consists of generating vehicle routes to serve a set of prespecified requests while minimizing the total routing cost during a certain planning horizon.

\quad D. Brevet et al \cite{bib4} dealt with the problem of on-demand transport with private vehicles and alternative locations (DARP-PV-AN), an extension that allows customers to manage their own trips using private vehicles. They add several alternative points for the origin/destination of transport requests to address privacy issues. A compact MILP model and an evolutionary local search algorithm (ELS) are proposed to solve the problem.

\quad Masson et al \cite{bib5} investigate the incorporation of transfers into DARP (DARPT). The objective is to minimise total routing cost under a set of additional transfer point-related constraints. They designed a formulation for the general case where the number of transfer locations is limited to one of each demand.

\quad Fleet homogeneity is a classic DARP assumption. It can be extended to heterogeneous fleets when customers have specific service requirements and expectations. The multi-depot variant with heterogeneous vehicles (MDHDARP) was presented by Detti et al. \cite{bib6}, the studied problem derives from a real-world healthcare application. Among the constraints, we can list vehicle capacity, pickup and delivery time windows, precedence constraints, and quality of service provided. They extended the problem by introducing vehicle-patient compatibility to ensure that the assigned vehicle is compatible with the patient's needs. In addition to heterogeneous fleet management, Parragh et al. \cite{bib21} consider heterogeneous users and constraints related to the drivers, such as drivers’ work duration limits and lunch breaks. To tackle this problem, they implemented a hybrid algorithm combining column generation with a VNS heuristic algorithm.

\quad DARPs belong to the family of complex combinatorial problems whereas concerning most often many requests. To address those problems, exact methods find difficulties that drive several contributions to the use of metaheuristics. Combining more than one metaheuristic results in hybrid metaheuristics that have been proven to perform better than the classic ones. In the following, we concentrate on the meta-heuristics used in the design of our approach. For each metaheuristic, we concentrate on important contributions.

\quad The Large Neighborhood Search (LNS), Masson et al. \cite{bib5} used an adaptive strategy into LNS (ALNS) with destroy and repair operators to deal with (DARPT). The adaptiveness of the heuristic is ensured by a mechanism that calls each operator depending on its performance in previous iterations. Molenbruch et al. \cite{bib7} designed a large neighborhood search metaheuristic to solve the collaborative DARP.

\quad Simulated Annealing (SA), Braekers et al. \cite{bib8} presented a version of Simulated Annealing to solve large instances for the MD-H-DARP (noted as Deterministic Annealing). They proved that the algorithm able to produce good results, both in terms of solution quality and computation time. Reinhardt et al.\cite{bib9} proposed a simulated annealing-based heuristic to deal with a DARP for disabled passengers at airports with complicating synchronization constraints.

\quad The Iterated Local Search (ILS), unlike the above-mentioned metaheuristic, the iterated local search is not widely used to solve the DARP. The only work that we found was that of Malheiros et al.\cite{bib10}. The algorithm proposed combines the iterated local search algorithm with a set partitioning approach (ILS-SP) and obtains high-quality solutions of multi-depot heterogeneous DARP within a short time.

\quad Path-relinking (PR) is a common ingredient of hybrid optimization methods. Parragh et al. \cite{bib11} used it as a post-optimization within an iterated variable neighborhood search algorithm to tackle a bi-objective DARP. In Molenbruch et al. \cite{bib7} a path relinking phase was applied after a multi-directional local search algorithm to deal with a DARP with real-life characteristics concerning the medical needs of patients.

\begin{sidewaystable}
    \scriptsize
	\centering
	\caption{Literature review summary}\label{tab_review}
	\begin{tabular}{|l|l|l|l|}
		\hline Reference & Objectif Function & Constraints & Solution Methods\\
		
		\hline 	
		Parragh et al.\cite{bib2}
		& \begin{tabular}{l} 
			Min Routing cost \\
		\end{tabular} 
	    & \begin{tabular}{l} 
	     Standard DARP constraints  \\
	    \end{tabular}
		& \begin{tabular}{l} 
			Variable neighborhood search (VNS) \\
			\hspace{1.5cm} + \\
			Path relinking (PR)\\
		\end{tabular} \\
		\hline
		
		Parragh et al.\cite{bib21}
		& \begin{tabular}{l} 
			Min Routing cost \\
		\end{tabular} 
		& \begin{tabular}{l} 
			Standard DARP constraints\\
			\hspace{1.5cm} + \\
			Heterogeneous clients/vehicles\\
			\hspace{1.5cm} + \\
			Max route duration\\
			\hspace{1.5cm} + \\
			Lunch breaks
		\end{tabular}
		& \begin{tabular}{l} 
			Variable neighborhood search (VNS) \\
			\hspace{1.5cm} + \\
			Column generation algorithm (CG)\\
		\end{tabular} \\
		\hline
		
		Reinhardt et al.\cite{bib9}
		& \begin{tabular}{l} 
			Min Routing cost \\
			\hspace{1.5cm} + \\
			Max Number of met requests \\
		\end{tabular}
		& \begin{tabular}{l} 
			Standard DARP constraints\\
			\hspace{1.5cm} + \\
			Transfers\\
		\end{tabular}
	     & \begin{tabular}{l} 
	         Simulated Annealing (SA)  \\
	     \end{tabular} \\ 
		\hline 
		
		Masson et al.\cite{bib5}
		& \begin{tabular}{l} 
			Min Routing cost\\
		\end{tabular} 
		& \begin{tabular}{l} 
			Standard DARP constraints\\
			\hspace{1.5cm} + \\
			Transfers\\
		\end{tabular} 
	    & \begin{tabular}{l} 
	    	Adaptive Large Neighborhood Search (ALNS)\\
	    \end{tabular} \\ 
		\hline 
		
		Braekers et al.\cite{bib8}
		& \begin{tabular}{l} 
	      	Min Routing cost\\
		\end{tabular}
		& \begin{tabular}{l} 
			Standard DARP constraints\\
			\hspace{1.5cm} + \\
			Multi-depot\\
		\end{tabular} 
		& \begin{tabular}{l} 
			Branch and Cut\\
			\hspace{1.5cm} + \\
			Deterministic Annealing (DA)\\
		\end{tabular}\\
		
		\hline
		Molenbruch et al.\cite{bib7}
		& \begin{tabular}{l} 
			Min Routing cost\\
			\hspace{1.5cm} + \\
			Min Mean user ride times\\
		\end{tabular}
		& \begin{tabular}{l} 
			Standard DARP constraints\\
			\hspace{1.5cm} + \\
			Restrictions user–user/user–driver\\
		\end{tabular} 
		& \begin{tabular}{l} 
			Path Relinking (PR)\\
			\hspace{1.5cm} + \\
			Multi-directional local search algorithm\\
		\end{tabular} \\
		
		\hline
		Malheiros et al.\cite{bib10} 
		& \begin{tabular}{l} 
			Min Routing cost\\
		\end{tabular} 
		& \begin{tabular}{l} 
			Standard DARP constraints\\
			\hspace{1.5cm} + \\
			Multi-depot\\
			\hspace{1.5cm} + \\
			Heterogeneous vehicles\\
		\end{tabular} 
		& \begin{tabular}{l} 
			Iterated Local Search (ILS)\\
			\hspace{1.5cm} + \\
			Set partitioning (SP)\\
		\end{tabular}\\
		
		\hline	
		
		Parragh et al.\cite{bib11} 
		& \begin{tabular}{l} 
			Max Total profit\\
		\end{tabular} 
		& \begin{tabular}{l} 
			Standard DARP constraints\\
			\hspace{1.5cm} + \\
			Split loads constraints\\
		\end{tabular} 
	    & \begin{tabular}{l} 
	    	Branch and Price (B\&P)\\
	    \end{tabular} \\
		\hline
		Masmoudi et al.\cite{bib3} 
		& \begin{tabular}{l} 
			Min Routing cost\\
			\hspace{1.5cm} + \\
			Min BSS costs\\
		\end{tabular}
		& \begin{tabular}{l} 
			Standard DARP constraints\\
			\hspace{1.5cm} + \\
			G-VRP constraints\\
		\end{tabular}
	   & \begin{tabular}{l} 
	   	Evolutionary Variable Neighborhood Search (EVO-VNS)\\
	   \end{tabular} \\
		\hline
		Brevet et al.\cite{bib4}
		& \begin{tabular}{l} 
			Min Routing cost\\
		\end{tabular}  
		& \begin{tabular}{l} 
			Standard DARP constraints\\
			\hspace{1.5cm} + \\
			Private vehicles\\
			\hspace{1.5cm} + \\
			Alternative nodes\\
		\end{tabular}  
	    & \begin{tabular}{l} 
	    	Evolutionary Local Search (ELS)\\
	    \end{tabular} \\    
		\hline
		
		Detti et al.\cite{bib6} 
		& \begin{tabular}{l}
			Min Total distance traveled \\
			\hspace{1.5cm} + \\
			Min Total travel time \\
			\hspace{1.5cm} + \\
			Min vehicles waiting times \\
		\end{tabular} 
		& \begin{tabular}{l}
			Standard DARP constraints\\
			\hspace{1.5cm} + \\
			Multi-depot\\
			\hspace{1.5cm} + \\
			Heterogeneous vehicles\\
			\hspace{1.5cm} + \\
			Compatibility constraints\\
		\end{tabular} 
		& \begin{tabular}{l}
			Variable Neighborhood Search (VNS)\\
			\hspace{1.5cm} + \\
			Tabu Search (TS)\\
		\end{tabular} \\  
		\hline
	\end{tabular}
\end{sidewaystable}

\quad Despite the large variety of DARP variants arising from real-life scenarios, to the best of our knowledge, there is still some gap in driver requirements when designing schedules. In fact, there is a lack of papers that address constraints related to drivers, which are typically defined in terms of work duration limits and lunch breaks. With respect to the literature on DARP variants, we propose to take driver availability and preferences into account while scheduling requests. Thus, as distinguish from the standard DARP formulation, we incorporate a threshold constraint and a lateness penalty within the objective function to ensure the accurate modeling of the driver’s arrival time to a specific area of interest. In terms of solution methods, there have been successful hybrid approaches proposed for different DARP variants. However, most of these methods involves only two methods, which motivated us to design a hybrid algorithm that thoroughly explores the benefits of several approaches. 

\quad The table \ref{tab_review} provides an overview of the literature review with the column Objectif Function presenting the criteria incorporated in the model's Objectif Function. The column Constraints displays the constraints considered in the model (the standard constraints refer to the common constraints used in the model presented by Cordeau and Laporte (2003)). The column Algorithms describes the articles' methodologies.

\section{DARPDP: Description and Formulation}\label{Sec3}
	
\quad We propose an extension of the Dial-A-Ride problem that takes advantage of shared taxis and carpooling. Typically, we consider a set of clients who want to share a taxi and geographically dispersed taxis. Given that each driver's taxi has preferences, we assume that starting from a location representing his current position, he seeks to achieve a profitable route in an area close to his final destination. This destination is considered as known high demand at a specific time,  such as (bus stops, streetcar stations, electric buses, and train stations). From the driver's perspective, we try to model the problem so that he can benefit from every time saving while achieving his objectives and without ignoring the clients' preferences. A solution for the problem is represented by a set of routes. We base the assumptions of the model on three main components: client, driver, and route:
 
\textbf{Client (request)} we associate each request with two vertices a vertex of origin and a vertex of destination, and time windows associated with these two points. Additionally, a demand representing the number of passengers to transport as well as a service time, and a maximum ride time constrains the duration of his trip.

\textbf{Driver (taxi)} To each taxi, we associate two vertices, an origin vertex and a destination vertex with a time window.
     
\textbf{Route} Each route respects many constraints: the capacity of the taxi, a restricted duration, the order of passage on the vertices of each request, and the time windows associated with this later (hard constraints). For the taxi destination time window, we assume it is hard with a small epsilon tolerance. For example, as commonly found in train stations, even if the taxi driver misses the arrival of a train, the crowd leaving the station does not disperse on the spot, and a small tolerance has no repercussions.

\begin{figure}[h!]
	\includegraphics[width=0.96\textwidth,height=0.8\textwidth]{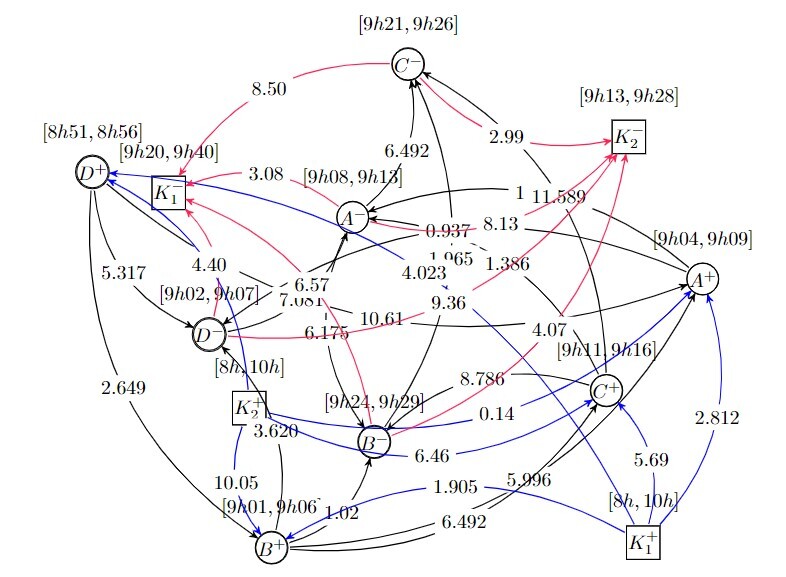}
    \caption{Illustrative example of the DARPDP network}\label{fig1}
\end{figure}

\quad Since this problem is an extension of the DARP, we adopt a model previously defined in the literature (Cordeau's formulation \cite{bib1}), where we adjust the objective function and the constraints to improve driver satisfaction. To describe the mathematical model we use the notations listed in Table \ref{tab_notations}.
\begin{table}[!ht]
   \begin{center}
      \caption{List of notations}\label{tab_notations}%
    	\begin{tabular}{ll}
    	\hline
    	   Notation  & Description\\
    	\hline
            \underline{Sets}\\
            $\mathcal{G}=(\mathcal{V}, \mathcal{A})$   &  An oriented and weighted graph illustrates the given network.\\
            $\mathcal{A}$ & Set of arcs.\\
            $\mathcal{V}$ & Set of all vertices $\mathcal{V} = \{s_{k}/k \in \mathcal{K}\} \cup \mathcal{P} \cup \mathcal{D} \cup \{t_{k}/k \in \mathcal{K}\}$.\\
             $\mathcal{N}$ & Set of all requests' vertices $\mathcal{N} = \mathcal{P} \cup \mathcal{D}$.\\
            $\mathcal{K}$ & Set of all taxis $\mathcal{K} := \{1, . . ., m\}$.\\
            $\{s_{k}/k \in \mathcal{K}\}$ & Set of taxis' origins $\{s_{k}/k \in \mathcal{K}\} := \{1, . . ., m\}$.\\
		$\{t_{k}/k \in \mathcal{K}\}$ & Set of taxis' destinations $\{t_{k}/k \in \mathcal{K}\} := \{n + m, . . ., 2m + n\}$.\\
            $\mathcal{P}$ & Set of pickups' vertices $\mathcal{P} := \{m, . . ., m +n\}$\\
		$\mathcal{D}$ & Set of deliveries' vertices $\mathcal{D} := \{2m+n, . . ., 2(n + m)\}$\\

            \underline{Parameters}\\
    				  	
    		$n, m$  &   The total number of requests, taxis respectively.\\
			$ s_{i}$ & A service time of vertex i $( s_{i} = 0 \quad, \forall i \in \{s_{k} / k \in 
                \mathcal{K}\} \cup \{t_{k} / k \in \mathcal{K}\})$.\\
			$[e_{i}, l_{i}]$ & Represent the earliest and the last time at which the service can begin\\ 
			& at vertex i.\\
			$ q_{i}$ & A demand of client i (number of passengers to be picked 
                            up or delivered)\\
                & $(q_{i} = 0, \forall i \in \{s_{k} / k \in \mathcal{K}\} \cup \{t_{k} / k \in \mathcal{K}\}$ and $q_{i} = -q_{i+m}, \forall i \in \mathcal{P})$ \\
			$d_{i,j}$ & The distance required to cross the arc $(i, j)$.\\
		    $t_{i,j}$ & The time required to cross the arc $(i, j)$.\\
		    $T_{r}$ & The maximum duration of the routes of any taxi.\\
		    $T_{M}$ & The maximum ride time of each request.\\

                \underline{Decision variables}\\        
                $x_{i,j}^{k}$ & Binary variable equal to 1 if vertex j is served immediately after vertex i\\
                & by a taxi k, and to 0 otherwise;\\
		        $T_{i}^{k}$ & Continuous variable corresponding to the time at which taxi k starts\\
                & its service at vertex i;\\
		    $Q_{i}^{k}$ & Continuous variable corresponding to the number of passengers in the taxi k\\
                & after servicing location i of $\mathrm{V}$;\\
		    $y_{i}^{k}$ & Binary variable, equal to 1 if client i is served by a vehicle k (each client\\
                & is represented by his pickup point i);\\                   
    		\hline
       \end{tabular}
     \end{center}
\end{table}

Therefore, we base the objective function on the following three criteria:
\begin{enumerate}
\item Maximise the number of served requests, which is equivalent to minimizing the total number of rejected requests; we intend to satisfy as many requests as possible during the duration of the route.
\item Minimise the total distance traveled (expressed in miles).
\item Set a tolerance value of  $\epsilon_{k}$ for the driver's destination time window and minimise any deviation from this latter (expressed in minutes).
\end{enumerate}
\quad There are several algorithms for solving multi-objective optimization problems. We consider the hierarchical and the weighted sum methods. We use the first one to express the priority of maximizing the number of served requests so that the solution doesn't sacrifice more than $20\%$ of its optimal value to optimize the other objectives (which means a solution to the problem must satisfy at least $80\%$ of the total requests). To tackle the remaining objectives, we use a weighted sum approach.    
\quad Indeed, it is necessary to make sure that all terms in the objective function are comparable and have the same unit of measurement. One suggestion is to express the objective function in terms of time. In the first stage, we convert the distance from miles to kilometers and then to minutes by multiplying it by the speed (One note is that while the distance is now expressed in minutes, it doesn't accurately reflect the total duration of the route because it doesn't take into consideration the waiting times at different stops.\\
\quad The mathematical model is represented as follows:\\
\emph{The Objective-Function}
\begin{equation}
 Min (w_{1}\times\sum_{k \in K} \sum_{i\in V} \sum_{j \in V}  c_{i,j}^{k} x_{i,j}^{k}+ w_{2}\times\sum_{k \in K}max\{0, T_{t_{k}}^{k}-l_{t_{k}}\})\label{eq1} 
\end{equation}
subject to the constraints:
\begin{align}
	&\sum_{i \in \mathcal{P}}(1 - \sum_{k \in \mathcal{K}} y_{i}^{k}) \leq \epsilon \label{eq}\\
	&\sum_{k \in \mathcal{K}} y_{i}^{k} \leq 1 && \forall i \in \mathcal{P}\label{eq}\\
	&\sum_{k \in \mathcal{K}} \sum_{j \in \mathcal{V}} x_{i,j}^{k} = \sum_{k \in 
		\mathcal{K}} y_{i}^{k} && \forall i \in \mathcal{P}\label{eq}\\
	&\sum_{j \in \mathcal{V} \backslash{t_{k}}} x_{j,i}^{k} - \sum_{j \in \mathcal{V} \backslash{t_{k}}} x_{j,i+m}^{k} = 0 && \forall i \in \mathcal{P},  \forall k \in \mathcal{K}\label{eq}\\
	&\sum_{j \in \mathcal{V} \backslash{\mathcal{D}}} x_{s_{k},j}^{k}=1 && \forall k \in \mathcal{K}\label{eq}\\
	&\sum_{j \in \mathcal{V}\backslash{t_{k}}} x_{j i}^{k}=\sum_{j \in \mathcal{V}\backslash{s_{k}}} x_{i j}^{k} && \forall i \in \mathcal{N}, k \in \mathcal{K}\label{eq}\\
	&\sum_{i \in \mathcal{V}\backslash{\mathcal{P}}} x_{i,t_{k}}^{k}=1 && \forall k \in \mathcal{K}\label{eq}\\
	&T_{i}^{k}+s_{i}+t_{i, j} \leq T_{j}^{k}+M\left(1-x_{i j}^{k}\right) && \forall i \in \mathcal{V},j \in \mathcal{V}, k \in \mathcal{K}\label{eq}\\
	&t_{t, i+m} \leq T_{i+m}^{k}-\left(T_{i}^{k}+s_{i}\right) \leq T_{M} && \forall i \in \mathcal{P}, \forall k \in \mathcal{K}\label{eq}\\
	&e_{i} \leq T_{i}^{k} \leq l_{i} && \forall i \in \mathcal{N}, k \in \mathcal{K}\label{eq}\\
\end{align}
\begin{align}
	&T_{t_{k}}^{k} \leq l_{t_{k}} + \epsilon_{k} && \forall k \in \mathcal{K}\label{eq}\\
	&Q_{i}^{k}+q_{j} \leq Q_{j}^{k}+M\left(1-x_{i j}^{k}\right) && \forall i \in \mathcal{V},j \in  \mathcal{V}, k \in \mathcal{K}\label{eq}\\
	&max\{0, q_{i}\} \leq Q_{i}^{k} \leq min\{Q_{k}, Q_{k}+q_{i}\} && \forall i \in \mathcal{V},k \in \mathcal{K}\label{eq}\\
	&Q_{s k}^{k}=Q_{t_{k}}^{k}=0 && \forall k \in \mathcal{K}\label{eq}\\
	&x_{i j}^{k} \in\{0,1\} && \forall i \in \mathcal{V}, j \in \mathcal{V}, k \in \mathcal{K} \label{eq}
\end{align}
            
\paragraph{} The objective function (1) aims to minimize the total sum of 
the distances traveled and the violation of the driver’s desired destination time. Constraint (2) is the $\epsilon-constraint$ that bound the total number of unserved requests. Constraints (3) ensure that each request is served at most once by a unique taxi. Constraints (4) ensure consistency between the variables $x_{i,j}^{k}$ and $y_{i}^{k}$. Constraints(5) guarantee that the request's pickup and delivery vertices are served by the same vehicle. Constraints (6\_8) guarantee that each taxi starts its trip at its current position, leaves each vertex it arrives at, and ends at the destination point that corresponds to it. The consistency of the time is ensured through constraints (9). Constraint (10) ensures that the ride time for each request must not exceed the maximum ride time allowed. Constraints (11) guarantee that all requests are handled within their time windows. Constraint (12) ensures the driver’s preferences since the time window of the destination is assumed to be soft but controlled. We enable an increase by an epsilon that represents the minimum value by which we can exceed the driver’s desired time. Otherwise, it will be penalized in the objective function. At each vertex, constraints (13\_15) were designed to make sure that the sum of passengers in the taxi is smaller than its capacity. Constraint (16) states the binary domain.

 \section{The proposed Method}\label{Sec4}
	
\quad The proposed method to solve the DARPDP is a metaheuristic method based on an Iterated Local Search ILS scheme which was introduced for the first time by Lourenço et al. (2003) \cite{bib12}. ILS has significant success in solving vehicle routing problems in the literature. For a complete description of this method, we refer to (Lourenço. Martin and Stützle (2010) \cite{bib13}). The following key points constitute the basis of the proposed approach strategy : 
	
\begin{itemize}
    \item A two phases heuristic to generate initial solutions;
    \item A local search to improve the current solution;
    \item A post-optimization method to efficiently traverse the solution space;
    \item An acceptance Criterion to define whether a new solution should be accepted;
    \item A perturbation to diversify the search space.
\end{itemize}

We give the general framework's pseudo-code in Algorithm  \ref{alg: E-ILS} and further details on each component in the following sections.

\quad Algorithm \ref{alg: E-ILS} gives the general structure of the proposed method, which starts with the initialization of different parameters and the construction of the initial solution. The main loop in lines (6-26) performs until reaching a maximum run time. At each iteration, a learning-based local search phase improves the current solution. The path relinking then enhances the solution. If the resulting solution $s"$ is better than the current one $s$, $s"$ replaces $s$, and if it is better than the best, $s"$ replaces $s^{*}$. Otherwise, the number of iterations without improvement increments, and a call to simulated annealing is required to determine whether to accept $s"$ as a new solution. When an even number of iterations without improvement is met, the algorithm executes a Perturbation phase. Finally, The temperature parameter T decreases and restarts after reaching a minimum temperature.

\begin{algorithm}[h!]
    \caption{The proposed approach E-ILS}\label{algo1}
	\begin{algorithmic}[1]
	      \Require
			\Statex Input Data : Contains information about vertices, requests and vehicles
			\Statex $CPU_{max}$: The maximum running time
			\Statex $T_{min}, T_{max}$:The minimum and the maximum temperature, respectively
			\Statex $sizeN$ : The size of the neighborhoods
			\Statex $sizeE$: The maximum size of the set E, the set of elite solutions    
			\Ensure $s^{*}$	 
			
			\State $s_{0} \leftarrow TwoPhasesHeuristic(Input Data)$
			\State $s, s^{*} \leftarrow s_{0}$
			\State $T \leftarrow T_{max}$
			\State $NoImprov  \leftarrow 0$
			\State $Update\_setE(s, sizeE)$
			
			\While{$ CPU < CPU_{max}$}
			
			     \State $s' \leftarrow Learning Based Local Search(s, s^{*}, sizeN)$
			     
			     \State $s'' \leftarrow PathRelinking(s', E)$	
			
			     \If{$F(s'')<F(s)$}
			           \State$ s \leftarrow s''$
			           \If{$F(s'')<F(s^{*})$}
			                 \State$ s^{*} \leftarrow s''$
			            \EndIf
			     \Else
			           \State $NoImprov \leftarrow NoImprov + 1$
			           \State $s \leftarrow SimulatedAnnealing(s, s'', T)$
			      \EndIf
			      
			      \State $Update\_setE(s, sizeE)$
			      
			      \If{$NoImprov \bmod 2$} 
		     	        \State$ s \leftarrow Perturbation(s, sizeN)$
			       \EndIf
			       
			       \State $T \leftarrow T*\alpha_{T}$
			           \If{$T <= T_{min}$}
			                \State$ T \leftarrow T_{max}$;
			           \EndIf
		   	      \EndWhile      
    \end{algorithmic}\label{alg: E-ILS}
\end{algorithm} 

\subsection{The initial solution}\label{subsec_1}
	
\quad In our proposed approach, we generate the initial solution by a new two-phase heuristic based on a cluster-first route second principle. First, we perform a "clustering" phase, after which the requests are decomposed based on spatial characteristics. The idea is that requests will be assigned to the taxi if they are located in a specific area (assumed around the direction in which the driver wants to drive). To this end, we approximate the route taken by the taxi by an imaginary line connecting its origin and destination points. After assigning requests to routes, we perform sequential and parallel insertions to insert requests into routes. We explain the method in detail in the following.

\begin{algorithm}[!ht]
    \caption{TwoPhasesHeuristic}
	   \begin{algorithmic}[1]	
	       \Require Input Data: V(vertices), R(requests), K(taxis)
	       \Ensure  An initial solution s
	    		
	       \For{$(v, k) \hspace{3pt}\in (V,K)$}
	    		\State Compute $d_{v,k}$ 
	       \EndFor
	    		
	       \For{$k\in K$}
	    		\State Compute $\bar{d_{k}} = \frac{\sum_{v \in V}^{} d_{v,k}}{|V|}$ 
	    		\State $L_{k} = \{ r\in R/ d_{r+,k} < \bar{d_{k}}, d_{r-,k}< \bar{d_{k}}\}$
	       \EndFor
              \For{$r\in R$}
	    		\If{$r\in \bigcap_{i \geq 2}L_{k_{i}}$}
	    		   \State $I_{r} = \{ i / r\in \bigcap_{i}L_{k_{i}}\}$
	    		     \State $\{ L_{k_{i}} / k_{i} \in K, i\in I_{r} \} = \{ L_{k_{i}}/\{r\} / k_{i} \in K, i\in I_{r} \}$
	    		   \State $L_{MC} \leftarrow L_{MC} \cup \{ (r, L_{rk}) / rk \in K\}$
	    		\EndIf
	       \EndFor
	       \State $ s_{0} \leftarrow An empty solution $
	       
           \Statex \textbf{\emph{Sequential insertion}}
           
	       \For{$route\in s_{0}$}
	    		\State Sort $L_{k}$ by the upper bound for the start of service of each request 
	    		\State ($min\{l_{p},l_{d}-T_{p,d}- s_{d}\}$)        
	    		\For{request r in $L_{k}$}
	    		  \State $route_{k} \leftarrow insert\_best(r, route_{k})$.
	    		\EndFor   
	        \EndFor 
	    		    
	       \Statex \textbf{\emph{Parallel insertion}}
	    	\For{request rc in $L_{MC}$}
	    		\State Choose the most profitable route.
	    		\State $best\_route_{rc} \leftarrow insert(rc, best\_route_{rc})$
	       \EndFor		    
	       \State For the step 6 and 7 if an insertion is failed the request will be rejected.
	    	\end{algorithmic}\label{alg: TPH}
	    \end{algorithm}
	    
\quad The algorithm \ref{alg: TPH} starts by calculating the distances between the vertices and different imaginary lines (lines 1-3). We use the minimum distance from a point to the line segment as a distance measure and the euclidean metric. The average distances to each taxi's line are then used to create a candidate list of the requests that can be assigned to this taxi (lines 5-6). A request will be assigned to a taxi if its origin and destination are within a certain radius around the taxi direction. One note is that a request can be part of multiple clusters (we associate each cluster with a list). For this, we differentiate between a list specific to a taxi and one containing the requests belonging to multiple clusters (we note the list $L_{MC}$). In the latter case (lines 8-14), we save the index of these clusters in a list noted as $I_{r}$ specific to each request. The request will then be removed from its original list (cluster) and inserted into the list $L_{MC}$ that contains pairs of (request, potential route assignment). In lines (16-26), the routes' construction phase is carried out, starting with the initialization of an empty solution (each route starts at the driver's origin location and ends at his destination location). We proceed to insert sequentially, in each route, its specific candidate requests list in the best position. Therefore, we insert in parallel the remaining requests in $L_{MC}$ in different potential routes and choose the best route. Note that the best possible insertion is the one that satisfies all the constraints and minimizes the distance the most.

\subsection{Learning local Search}\label{subsec_2}
	
\quad The purpose behind this local search technique is to enhance a given initial solution by focusing the search on promising areas. This method was inspired by the heuristic adopted in the work of Bo Peng et al. \cite{bib14}, which uses a learning process to determine which operators contribute most to the improvement of the solution. To do this, this mechanism applies a reward and penalty strategy to manage the movements dynamically. Thus, the operators that provide better results are chosen more often in the following iterations.
	
\textbf{Learning mechanism} At the beginning of the search phase, All moves are assigned the same score $(sc_{0} = 1)$ and thus the same probability of being chosen. After each iteration, the operator receives either a reward or a penalty in the form of scores, depending on the current solution and the best solution found so far. Then the new score is recorded, and the probability of choosing this operator is recalculated at the next iteration using the scores collected at the previous iterations. The probability is calculated by the formula:
\begin{equation}
      \lambda_{k}=\frac{sc_{k}}{\sum_{i=1}^{n} sc_{i}}, \quad k=1,2,\ldots, n.\label{eqLM}
\end{equation}
With n the total number of moves.\\
The operators used were selected after some experiments considering the impact on the final solution. They are briefly described as follows.
 
The \textbf{Relocate} consists of moving a request from one route to another (Braekers et al.\cite{bib8}).
	
The \emph{\textbf{Exchange}} consists in exchanging two requests from different routes. The best exchange is found by trying to insert the origin and destination of each request into the best position in the other route (Braekers et al.\cite{bib8}).
	
The \emph{\textbf{Exchange Natural}} consists of choosing two routes and performing all possible exchanges of a natural sequence in one route with any other natural sequence in another route. A natural sequence is a sequence of requests before and after which the vehicle is empty (Parragh et al.\cite{bib15}, Molenbruch et al. \cite{bib7}).
	
The \emph{\textbf{R-4-opt}} consists in selecting a route, then selecting four consecutive arcs in this route. These arcs are deleted, then the segments are reconnected by other arcs in a different order (Chassaing et al.\cite{bib16}).
	
The \emph{\textbf{Shift}} consists in choosing a route randomly. For all the combinations of two stops, and performs the exchanges that respect the precedence constraints.
	
At each call to an operator, the process behind it is repeated at most n times, where n is the size of the neighborhood. However, inserting all possible solutions would be very time-consuming. Thus, the algorithm stops the search after obtaining the first improving solution.

\begin{algorithm}[!ht]
    \caption{Learning Local Search}
        \begin{algorithmic}[1]
		\Require 
                 \Statex MaxIter: The maximum number of iterations without improvement
			     \Statex s: A solution to be improved
                 \Statex  sizeN
		\Ensure $s^{*}$		
		\State  $ I_{NoImprove} \leftarrow 0$
		\State  $ \mathcal{N} \leftarrow$ The set of neighborhood structures to be used. 
		\State  $ Scores \leftarrow \{ score_{N} = 1 / N\in \mathcal{N}\}$
  
		\While{$I_{NoImprove} < MaxIter$}
		    \State $ P \leftarrow \{ P_{N} = \frac{score_{N}}{\sum score_{N}} / N\in \mathcal{N}\}$
		      \State $ L \leftarrow \{ N\in \mathcal{N} / P_{N} > rand\}$
		      \State $L_{0} \leftarrow L$
		      \While{$ L \ne \emptyset$}
		           \State $ M \leftarrow rand(L)$
		           \State  $ s' \leftarrow M(s)$
		           \If{$f(s') < f(s)$}
		                  \State $s \leftarrow s'$
			            \State $L \leftarrow L_{0}$
			            \State $I_{NoImprove} \leftarrow 0$
			            \If{$f(s') < f(s^{*})$} 
			                \State $score_{M} \leftarrow score_{M} + \alpha * \beta_{1}$
			                \State $s^{*} \leftarrow s'$    
			            \Else
			                \State $score_{M} \leftarrow score_{M} + \alpha * \beta_{2}$
			            \EndIf
			            \Else  	
			                \State $score_{M} \leftarrow score_{M} * \gamma$
			                  \State $L \leftarrow L \backslash M $
			                  \State $I_{NoImprove} \leftarrow I_{NoImprove} + 1$     
			         \EndIf	
			      \EndWhile			
		       \EndWhile	
		\end{algorithmic}\label{alg: LLS}
	\end{algorithm}	
	
\quad The algorithm \ref{alg: LLS} starts with the initialization of the number of iterations without improvement and the list of moves. Then, assigns the value 1 to each operator to start the search. The main loop is performed until the number of iterations without improvement is reached. In the first entry in the loop, all the moves have the same probability of being called. As the execution progresses, this probability changes according to the performance of each neighborhood. In line 7, to build the list of candidate neighborhoods, a probability is chosen between 0 and 1, and another loop starts as long as this last list is not empty. Then, the procedure executes move by move randomly, and the resulting neighboring solution is noted as $s'$.Depending on the quality of this latter, the score will change; if this solution improves the current solution, the current solution will be modified. The search will be restarted by restoring the list of candidates and the number of iterations without improvement. Otherwise, the movement will be removed from the list and the number of iterations without improvement will be incremented by 1.
	
We can summarize the changes in the score as follow:
	
\begin{itemize}
    \item New best solution: whenever the use of an operator results in a new best solution, its score is increased by $(\alpha * \beta_{1})$:	
     $score \leftarrow score + \alpha * \beta_{1}$	
    \item New improving solution: each time the use of an operator results in a new solution that improves the current one, its score is increased by $(\alpha * \beta_{2})$:
    $score \leftarrow score + \alpha * \beta_{2}$
    \item Each time the use of an operator results in a new solution whose cost is worse than the current solution, its score is multiplied by $\gamma$.
    $score \leftarrow score * \gamma$
\end{itemize}
 
\subsection{Path-relinking}\label{subsec_3}

\quad After conducting tests on small instances solved optimally, we noticed that some solutions differ from one request in the wrong position from the optimal one. One thought is that a new iteration of learning local search can deviate from the global optimal solution. To address this, we further enhance the approach with a path relinking (PR) procedure. This permits fast intensification and leads to the desired solution in less time than a call to Learning Local Search.

\quad The (PR) was first proposed by Glover \cite{bib18} and applied as a post-optimization based on the exploration of paths connecting the current solution and the best elite solutions previously produced during the search. The reasoning behind this idea is that a path from the current solution to a guiding one will result in new solutions that combine the best features of the two solutions. We have implemented a back-and-forward strategy which is explained in algorithm 4.

\begin{algorithm}[!ht]
	\caption{Path-Relinking}
	\begin{algorithmic}[1]
		\Procedure{Path-Relinking}{$s^{i},s^{j}$}
		\State $\Delta \leftarrow \{j = 1, ..., n : s_{j}^{i} \ne s_{j}^{g}\}$
		\State $s^{*} \leftarrow argmin\{f(s^{i}), f(s^{g})\}$
		\State $f^{*} \leftarrow min\{f(s^{i}), f(s^{g})\}$
		\State $s \leftarrow s^{i}$
		
		\While{$|\Delta| > 1$}
		\State $l^{*} \leftarrow argmin\{ f(s \oplus l) : l \in \Delta\}$
		\State $ s \leftarrow s \oplus l^{*}$
		\State $ \Delta \leftarrow \Delta \backslash \{l^{*}\}$	
		\If{$f(s) < f^{*}$}
		\State $ s^{*} \leftarrow s$
		\State $ f^{*} \leftarrow f(s)$
		\EndIf
		\State $ s^{i} \leftarrow s^{g}$
		\State $ s^{g} \leftarrow s^{i}$
		\EndWhile
		\EndProcedure
	\end{algorithmic}\label{alg: PR}
\end{algorithm}

\quad In each call to (PR) \ref{alg: PR}, the current solution is considered as initial solution, and a guiding solution is randomly selected from the elite set E. The procedure starts with computing the set $\Delta$ (which represents the symmetrical difference between the two solutions and the set of movements required) (line 2). In lines (3-4), the algorithm states the best solution between the two solutions (the starting one) and (the guiding one) as $s^{*}$ and its cost value as $f^{*}$. And the current solution receives $s^{i}$ (the initial solution). The search then enters a while loop, which constitutes the main body of the algorithm. At each iteration, it carries out the following operations: generates a solution path connecting $s^{i}$ and $s^{g}$; deletes the move from $\Delta$ and returns the best solution along this path as the current solution. If this solution has an objective value lower than the best one, the best solution and its cost are updated. Next, the algorithm alternates between the initial solution $s^{i}$ and the guiding solution $s^{g}$, and a new iteration starts. The while loop continues until $\Delta$ becomes empty.
	
\quad Then a set $\Delta$ composed of requests in the initial solution that are considered as being in the wrong route with respect to the guiding solution. The set $\Delta$ is computed based on these three key moves as in Santos and Xavier \cite{bib19}:	
$$\left\{
\begin{array}{l}
\text{1. } r \notin s_{i} \text{ and } r \in s_{g} \\
\text{2. } r \in s_{i} \text{ and } r \notin s_{g} \\
\text{3. } r \in s_{i} \text{ and } r \in s_{g} \text{ but } id_{R_{s_{i}}^{r}}\neq id_{R_{s_{g}}^{r}}
\end{array}
\Rightarrow
\left\{\begin{array}{l}
\text{1. Insert r in }s_{i}\slash id_{R_{s_{i}}^{r}}=id_{R_{s_{g}}^{r}}\\ 
\text{2. Remove r from }s_{i}\\
\text{3. Remove r from }R_{s_{i}}^{r} \text{and insert}\\
\text{it in }R_{s_{i}}^{\prime r}\slash id_{R_{s_{i}}^{\prime r}}=id_{R_{s_{g}}^{r}}\\ 
\end{array}
\right.
\right.$$

Where $r$ a request, $s_{i}$ the initial solution , $s_{g}$ the guiding solution, $id_{R_{s_{i}}^{r}}$ represents the id of the route $R$ associated with a request $r$ in the solution $s_{i}$.

\subsection{Perturbation}\label{subsect_4}

\quad During the process of generating the initial solution, we did not focus on reducing the number of rejected requests, and it turns out that the principle behind the local search algorithm results in an improved solution that maintains the same number of requests as the initial solution. More specifically, the operators can improve the current solution by making changes that affect only the requests that are part of the solution. However, when a request is removed from the current solution, we consider the possibility of inserting others. One approach could be to take advantage of the perturbation phase to delete less profitable requests and insert others. 

\quad Therefore, the perturbation can be seen as a Large Neighborhood Search (LNS) based approach. The principle of this approach is to remove a number of requests from their current routes and rearrange them at an optimal cost. We opt for this method in order to attempt to insert requests that were rejected in prior iterations. To do this, we begin by removing a certain number of requests from the current solution. Next, we reconstruct the solution by iteratively inserting as many of the previously rejected requests as possible, using repair operators.\\
We have implemented four removal operators (Route Removal, Greedy Removal, Random Removal, and Related Removal) and three repair operators (Compatibility Insertion, Best Insertion, and Random Insertion).
\subsubsection*{Destruction operators}
The objective of these operators is to delete a number of requests. We use four destruction operators: Route Removal, Greedy Removal, Random Removal, and Related Removal.
	
\textbf{\emph{Route removal}} The operator selects the route with the minimum number of requests and thus deletes all the requests figured in it (Nagata and Bräysy \cite{bib20}.
	
\textbf{\emph{Random removal}} The operator randomly selects n requests from the solution and then deletes them. All the routes modified will be optimized by the "Shift" operator (Pisinger and Ropke \cite{bib17}).
	
\textbf{\emph{Greedy removal}} For each request, the algorithm simulates its removal from its current route. The overall cost saving is calculated using a defined metric given in the following equation.

\begin{equation}
    sav(r) = 
    \begin{cases}
    (t_{p-1,p} + t_{p,d} + t_{d,d+1}) - t_{p-1,d+1} \quad if \quad d=p+1\\
    (t_{p-1,p} + t_{p,p+1} - t_{p-1,p+1}) + (t_{d-1,d} + t_{d,d+1} - t_{d-1,d+1}) \quad Otherwise\\
    \end{cases}
    \label{eqn: greedy_removal} 
\end{equation}
	
\quad Where p and d are the request's pickup and delivery positions, respectively. The $t_{i,j}$ represents the time spent in the solution to traverse the arc (i,j). The request whose deletion produces the largest saving is chosen. The operator repeats the same steps until n requests are deleted.
 
\textbf{\emph{Related removal}} This operator removes a number of related requests (similar requests) from the solution. The measure of similarity is borrowed from data mining. We consider that each request is defined by four attributes (features); We note A as the set of these attributes. A = \{position of pickup point, the position of the delivery point, ready time of pickup point, due time of the delivery point\}. The value of a request $ r_{i} $ for an attribute $a_{k}$ is represented by the similarity function $ f(r_{i}, a_{k})$. The similarity between two requests $ r_{i} $ and $ r_{j} $ is denoted $ Sim(r_{i}, r_{j})$ and is represented by the aggregation of the similarity of all attributes divided by the number of attributes.
	
The similarity value between requests $r_{i}$  and  $r_{j}$ with respect to an attribute $ a_{k} \in A$ is defined as:
	
\begin{equation}
sim_{a_{k}}(r_{i}, r_{j}) = 1-
\frac{|f(r_{i}, a_{k}) - f(r_{j}, a_{k})|}{\max_{r \in R}f(r, a_{k}) - \min_{r \in R}f(r, a_{k})}
\end{equation}

(In order to normalize the values, we set them within the range 0-1 by dividing by the largest value (1 – value/largest value)).\\
Two requests $r_{i}$ and $r_{j}$ are related if the sum of the attribute values divided by the number of attributes is small. The first removed request $r_{i}$ is selected randomly. Then, the r-1 most related requests to $r_{i}$ are removed and added to the list of rejected requests.
	
\subsubsection*{Insertion operators}

To repair the solution after the destruction phase, we adapted three insertion methods to insert the rejected requests, namely, Best Insertion, Random Insertion, and Compatibility Insertion.
	
A simple insertion technique for the dial-a-ride problem involves inserting a request per iteration in parallel on different routes and examining all feasible insertions of both the pickup and delivery points. Therefore, two options are possible; choosing the best feasible insertion \emph{\textbf{Best Insertion}} or a random one \emph{\textbf{Random Insertion}}.
	
\emph{\textbf{Compatibility Insertion}} The Compatibility Insertion aims to insert the rejected requests in the most compatible route based on a compatibility matrix. To construct the matrix, we first calculate the distance between each pair of requests. The distance measure used is defined as follows:
	
\begin{equation}
    \begin{split}
        CompMeasure(r_1, r_2) & = max(\lvert l_{r_1^{+}} - e_{r_2^{+}} - t_{r_1^{+} , r_2^{+}} \rvert, \lvert l_{r_2^{+}} - e_{r_1^{+}} - t_{r_1^{+}, r_2^{+}} \rvert) \\
              & + max(\lvert l_{r_1^{+}} - e_{r_2^{-}} - t_{r_1^{+} , r_2^{-}} \rvert, \lvert l_{r_2^{-}} - e_{r_1^{+}} - t_{r_1^{+}, r_2^{-}} \rvert)\\
              & + max(\lvert l_{r_1^{-}} - e_{r_2^{+}} - t_{r_1^{-} , r_2^{+}} \rvert, \lvert l_{r_2^{+}} - e_{r_1^{-}} - t_{r_1^{-}, r_2^{+}} \rvert)\\
				& + max(\lvert l_{r_1^{-}} - e_{r_2^{-}} - t_{r_1^{-} , r_2^{-}} \rvert, \lvert l_{r_2^{-}} - e_{r_1^{-}} - t_{r_1^{-}, r_2^{-}} \rvert)
    \end{split}
\end{equation}\label{eqCI}

To determine the most compatible route for each request, we compute the sum of the distances between this request and all route requests, divided by the number of requests.
	
\begin{algorithm}[!ht]
   \caption{Perturbation}
	\begin{algorithmic}[1]
		\Procedure{Perturbation}{$s$}
		    \State Choose random r from the list of rejected requests.
		    \State $O^{-}, O^{+} := $ A set of destruction (repair) operators
		    \State $d \in O^{-}$ and $r \in O^{+}$
		    \State $s^{'} \leftarrow Destroy(s^{'}, d)$
		    \State $s^{'} \leftarrow Repair(s^{'}, r, R^{r}_{s^{'}})$	
		    \State \Return s'
		\EndProcedure
	\end{algorithmic}\label{alg:LNS}
\end{algorithm}
	
\quad The Perturbation is described in Algorithm \ref{alg:LNS}. Given "s" a current solution obtained by the local search. The algorithm chooses a number r of requests to remove from s (line 2). And two operators destroy and repair are selected from the list of operators. Then it degrades the current solution by removing r requests from s and putting them in the list of rejected requests. Then, the reinsertion operator completes the partial solution obtained by inserting as many requests as possible. 
	
\subsection*{Acceptance Criterion}\label{subsect_5}
	
\quad The last part of the heuristic is devoted to taking control of the acceptance of new solutions (i.e., determining whether the solution produced is accepted as the new current solution for the next iteration). Along the lines of many researches, we employed an acceptance criterion borrowed from simulated annealing. The basic idea is that, initially, we only accept solutions whose cost is lower than the cost of the current solution. And during the solution process, to better avoid getting trapped in a local minimum, deteriorating solutions are accepted with a given probability.
	
\begin{algorithm}[!ht]
   \caption{Based SA Acceptance Criterion}
	\begin{algorithmic}[1]
		\Procedure{SA}{$ s,s', T $}  
			\If{$f(s) < f(s')$}
			      \State \Return s
			\Else
			      \State $r \leftarrow rand(0, 1)$
			      \If{$r < \exp{\frac{f(s)-f(s')}{T}}$}
			        \State \Return s
			      \Else
			        \State \Return s'
			      \EndIf
		    \EndIf
	   \EndProcedure
     \end{algorithmic}\label{alg: SA}
\end{algorithm}	 
	
\quad The integration of simulated annealing: At the beginning of the search, algorithm \ref{alg: E-ILS} sets the temperature to its upper bound $T_{max}$. Once a new solution is returned, the procedure SA was called (algorithm \ref{alg: SA}) to decide with which solution the approach will continue the search. Besides the current temperature, the procedure takes as input the current solution $s'$ and the new solution $s$ to compare them. If $s$ is better than $s'$, it will return $s$. Otherwise, a random number r between 0 and 1 is chosen, and if $r < exp(-f(s)-f(s'))$ then returns $s$ otherwise returns $s'$. After this iteration, the current temperature is reduced by multiplying it with a given cooling rate $ \alpha $. If the current temperature reaches its lower bound $T_{min}$ (close to 0) then is reset to the initial value $T_{max}$.

\section{Computational Results}\label{Sec5}

\quad To validate our work, we first present the generation process of new test instances from real-world data of taxis in New York City. We next carry out a set of experiments using both algorithms proposed; an exact algorithm that results in solving the mixed integer mathematical programming and the approach proposed in section \ref{Sec4}.

The results were performed on an Intel(R) Core i7, with 2.60 GHz and 32 GB RAM, and the Windows operating system. The algorithms were coded in Python 3.8, where the model is solved using IBM ILOG CPLEX Optimization Studio 20.

\subsection{Test instances}\label{subsect_6}
	
\quad NYC Open Data is one of the largest and most widely used open data platforms in the world. We used the green and yellow taxi trip data from the year 2021. For a quick overview of the data set, the dataset has 22793 rows and 18 columns. A more detailed description can be found on the NYC open data website \url{https://opendata.cityofnewyork.us}.

The new benchmark problems for the DARPDP are generated from the yellow taxi trip data for January 2021. We only extracted the columns “PULocationID”, “DOLocationID”, “tpep\_pickup\_datetime”, “tpep\_dropoff\_datetime”, and the “trip\_distance”.We have considered a time horizon of 2 hours, and therefore we have restricted the “data\_times” columns so that only the “date\_time” between “08h10” and “10h10” are kept (this period was chosen because of an extremely interesting demand).

We constructed 28 instances (where the number of requests in \{10, 15, 20, 30, 50, 100, 200\} according to 4 scenarios. The scenarios are generated according to the following characteristics: the width of the time windows of the customers (either small of 5 minutes (type-a) or medium of 10 minutes (type-b)); the number of taxis to be used, calculated as nbr/9 (with "nbr" the number of requests, and we assume that each vehicle will serve an average of 9 passengers in its trip); we took either the lower or upper bound.

\quad In the following, we describe the process for generating each instance from the original reference dataset:

\textit{Choice of the locations of the client vertices (Pickups/ Deliveries) and the vehicle vertices (Origins/ Destinations)} We joined the Lat/Long geographic coordinates on the "PULocationID" and "DOLocationID" columns by merging the dataset with the shape files for the taxi zones file. 

\textit{Choice of the clients' time windows} For setting the time windows for the different vertices, we used the column “tpep\_pickup\_datetime” and “tpep\_dropoff\_datetime” as the center of the time windows associated with the pickup and drop-off points with the consideration of the width (either small or medium).
    
\textit{Choice of the drivers' time windows} Whereas for the time windows of the taxis, for origins are flexible $ [490, 610] $, and for destinations, we generated the time windows by first choosing a uniform random number $ e_{i} $ in the interval $ [550, 595] $ in, then choosing $ (l_{i} = e_{i} + w)$ (with w the width of the time windows associated to the drivers' destination, w = 15 minutes).
    
\textit{Additional data} To make instances as realistic as possible, we consider that the maximum duration of each trip is 90 minutes and that the maximum time a passenger can stay in a taxi is 30 minutes. The service time is set at 30 seconds and represents the upper limit of time required for a passenger to get on and out of the taxi. And to simplify the problem, the demand is assumed to be 1, and the maximum capacity of all taxis is set to 3 passengers. All times are in minutes. The time between any two vertices of the network is determined using the Euclidean distance.

We refer to the instances as "inst\_an\_m" ("inst\_bn\_m"), representing type-a (type-b) instances with up to n requests and m vehicles. Regarding the number of requests, we divide the instances into three groups. Set small contains instances with (10 and 15 requests), the set medium has instances with (20, 30, and 50 requests), and set large contains instances with (100 and 200 requests). All test instances can be accessed from \url{https://github.com/SanaaOsd/TestInstancesDARPDP}.

Table \ref{table_instances} summarizes the main characteristics of the newly generated instances. The first column represents the instance name (either of type a or type b), the second column represents the number of requests to be served, and the third column the number of available taxis.

\begin{table}[!ht]
    \begin{center}
       \begin{minipage}{174pt}
        \caption{Instances and their characteristics}
    	\begin{tabular}{lcc}
    		\hline
    		Instance&Nb of requests&Nb of taxis\\
    		\hline		
    		inst\_10\_1   &      10 &    1\\	
    		inst\_10\_2   &      10 &    2\\
    		inst\_15\_1   &      15 &    1\\
    		inst\_15\_2   &      15 &    2\\
    		\hline
    		inst\_20\_2   &      20 &    2\\
    		inst\_20\_3   &      20 &    3\\
    		inst\_30\_3   &      30 &    3\\		
    		inst\_a30\_4  &      30 &    4\\
    		inst\_a50\_5  &      50 &    5\\
    		inst\_50\_6   &      50 &    6\\		
    		\hline
    		inst\_100\_11 &     100 &   11 \\
    		inst\_100\_12 &     100 &   12 \\
    		inst\_200\_22 &     200 &   22 \\
    		inst\_200\_23 &     200 &   23 \\		
    	    \hline
          \end{tabular}\label{table_instances}		
        \end{minipage}			
     \end{center}			
\end{table}

\subsection{Parameter settings}\label{subsect_7}

\quad The proposed method has several parameters to tune. We choose the parameter setting while developing the algorithm. Once the parameters related to a component(Method) are fixed to their best setting, we proceed to the next component.\\
\quad Starting with the learning local search, we choose to re-use the value of the size of the neighbourhoods ($sizeN$) used in the work of Parragh et al.\cite{bib15}, and the values reported in Bo. Peng et al. \cite{bib14} for the learning mechanism parameters, where ($\alpha$) is the reaction factor controls how quickly the score function reacts to changes according to the performance of the operators, ($\beta_{1}$) the reward parameter if an operator results in a new best solution, ($\beta_{2}$) the reward parameter if an operator improves the current solution, ($\gamma$) the penalty parameter if an operator results in a solution worse than the current one. Still the number of iterations ($Iter_{max}$), we investigated three settings where $Iter_{max}$ is chosen from (5, 10, 20).\\
\quad For parameters related to the SA method, we choose a commonly used value for the cooling rate ($\alpha_{T}$), and for the minimum and maximum temperature ($T_{min}$, $T_{max}$) we set candidate values depending on the number of requests and the number of vehicles n and m respectively, $\{(m/2, n/2), (m/2, n), (m, n)\}$.\\
Concerning the Path-Relinking method, it only requires a single parameter, the maximum size of the elite solutions, which is investigated with three settings (1, 5, 10). In order to obtain a perfect trade-off between the solution quality and the computation time, the maximum running time ($CPU_{max}$), is fixed regarding the number of requests in each instance, for the small type 300s, the medium instances (900s), and 1200s for the large ones. After some preliminary testing, the weights in the objective function are set to $w_{1} = 40$ and $w_{2} = 60$.\\
Table \ref{table_parameters1}, shows the results of 10 runs for each instance within a time limit of 300 seconds per run. For the $Iter_{max}$ parameter, we perform the learning local search only, for the SA parameter, we run the E-ILS without Path Relinking, and for the sizeE parameter, we run the E-ILS. The rows denoted "Best" and "Avg" represent the best and average solution values obtained with the algorithm, while the column "CPU" represents the average run time in seconds.\\
One can observe in the table \ref{table_parameters1} that setting the $Iter_{max}$ to 5 appear to require less CPU time, but with an increase in the solution value. In contrast, applying 20 may eventually lead to a slight improvement but at the expensive of more processing time, we decided to set $Iter_{max}$ to 10 since it offered a good balance between solution quality and CPU time. Setting the ($T_{min}$, $T_{max}$) to $(m/2, n/2)$ produces better solutions, in terms of both best and average results as well as computational time. For the $sizeE$ we conclude that 1 performs better than other values. Table \ref{table_parameters2} shows the retained values of the different parameters.

\begin{table}[!ht]
    \begin{center}
       \begin{minipage}{220pt}
        \caption{Settings of some important parameters for the E-ILS}\label{table_parameters1}
    	\begin{tabular}{llll}
    		\hline\noalign{\smallskip} 
    		$Iter_{max}$ & 5 & \textbf{10} & 20 \\
    		Best & 2757.6086 & 2727.1267 & 2725.0864 \\
    		Avg & 2839.62 & 2811.8222 & 2802.4756 \\
    		CPU[s] & 49.64567 & 51.70144 & 85.27037 \\
    		\hline\noalign{\smallskip}
    		$\left(T_{min}, T_{max}\right)$ & \textbf{(m/2, n/2)} & $(m / 2, n)$ & $(m, n)$ \\
    		Best & 1832.583 & 1848.905 & 1869.189 \\
    		Avg & 1951.851 & 1956.866 & 1964.608 \\
    		CPU[s] & 124.24 & 117.58 & 155.31 \\
    		\hline\noalign{\smallskip} 
    		$sizeE$ & \textbf{1} & 5 & 10 \\
    		Best & 1334.75 & 1388.041 & 1360.246 \\
    		Avg & 1433.261 & 1493.89 & 1445.335 \\
    		CPU[s] & 161.4597 & 175.5342 & 178.9134 \\
    		\hline
    	\end{tabular}		
        \end{minipage}			
     \end{center}			
\end{table}

\begin{table}[!ht]
	\begin{center}
		\begin{minipage}{174pt}
	\caption{Parameters setting adopted in the E-ILS}\label{table_parameters2}
		\begin{tabular}{ll}
			\hline\noalign{\smallskip}
			Parameter & Final Value\\
			\noalign{\smallskip}\hline\noalign{\smallskip}
			$\alpha$&  0.1 \\
			$\beta_{1}$& 5 \\
			$\beta_{2}$& 1\\
			$\gamma$& 0.9 \\                                
			$sizeN$&  3 \\
			$sizeE$&  1 \\
			$\alpha_{T}$& 0.99 \\
			$T_{min}$& $m/2$  \\	
			$T_{max}$& $n/2$  \\
			$CPU_{max}$& Small(300s)\\
			& Medium(900s)\\
			& Large(1200s)\\
			$Iter_{max}$& 10\\
			\noalign{\smallskip}\hline
		\end{tabular}
	\end{minipage}
\end{center}	
\end{table}

\subsection{Test Results}\label{subsect_8}

In this section, we describe the numerical experiments to evaluate the performance of the proposed method (noted as E-ILS) in solving DARPDP instances. We compare the proposed E-ILS algorithm to CPLEX for small instances and with three alternate versions of the method for the medium and the large instances. Besides comparing the objective function values and computational times of the obtained solutions, we evaluate the performance of the algorithms through three key performance indicators listed below:    

\begin{itemize}
\item The coverage rate represents the percentage of served clients.
\begin{equation}
      KPI_{1} = \frac{\sum_{k \in K}\sum_{i \in P} y_{i}^{k}}{\lvert P \rvert} \times 100\%
\end{equation}
\item The percentage of empty mileage refers to the distance traveled by vehicles without passengers onboard.
\begin{equation}
      KPI_{2} = \sum_{k \in K}\frac{\sum_{(i,j) \in V} T_{j}^{k} - T_{i}^{k}}{T_{t_{k}}^{k} - T_{s_{k}}^{k}} \times 100\%
\end{equation}  
Where i and j are two stops on which the vehicle drives empty $Q_{i}^{k} = 0$.

\item The excess ride time represents the average of the differences between the direct travel time between the pick-up location and the delivery location and the actual transport time of the client in the solution (in minutes).	
\begin{equation}
    KPI_{3} = \frac{ 
    \sum_{k \in K} \sum_{i \in P}T_{i+p}^{k}-(T_{i}^{k} + s_{i}) - t_{i, i+p}} {\lvert P \rvert}
\end{equation}   	
\end{itemize}

For each method, we solve each instance ten times and then compare the best objective function values and the computational times.

\subsubsection*{Comparison of CPLEX and the proposed E-ILS algorithm on small instances}

Table \ref{table_cplex_eils} illustrates the comparison between CPLEX and the proposed E-ILS algorithm for small DARPDP instances. The columns “Best” report the objective function value of the best solution obtained in 10 runs, "CPU[s]" refers to the computational time in seconds, and "gap[\%]" indicates the default gap as defined in CPLEX.

\begin{table}[!ht]
     	\begin{center}
     		\begin{minipage}{\textwidth}
     			\caption{Results of E-ILS and CPLEX on DARPDP instances}
     			\begin{tabular*}{\textwidth}{@{\extracolsep{\fill}}lccccc@{\extracolsep{\fill}}}
     				\hline
     				& \multicolumn{2}{@{}c}{CPLEX} & \multicolumn{3}{@{}c@{}}{E-ILS}\\
                        \cline{2-3}\cline{4-6}
     				Instance & {Best} & {CPU [s]} & {Best} & {CPU [s]} & {gap [\%]}\\
     				\hline
     				inst\_a10\_1   & 331.893 & 0.094  & 331.893 & 1.191   & 0 \\
     				inst\_a10\_2   & 137.1   & 0.204  & 137.1   & 16.145  & 0 \\		
     				inst\_a15\_1   & 272.922 & 3.719  & 272.922 & 27.969  & 0 \\		
     				inst\_a15\_2   & 318.586 & 0.281  & 318.586 & 17.807  & 0.003 \\		
     				inst\_a20\_2   & 410.026 & 0.531  & 410.026 & 93.501 & 0.01 \\
     				inst\_a20\_3   & 353.498 & 2.734  & 353.498 & 136.26 & 0.01 \\
         
                    inst\_b10\_1   & 172.758 &  3.719  & 172.758 & 4.47  & 0 \\
     				inst\_b10\_2   & 197.29  &  10.687 & 197.29  & 28.14 & 0 \\	
     				inst\_b15\_1   & 271.757 & 167.61  & 271.757 & 29.50 & 0.005 \\
     				inst\_b15\_2   & 333.858 & 72.75   & 333.858 & 120.82 & 0 \\	
     				inst\_b20\_2   & 253.495 & 584.141 & 253.495 & 349.72 & 0.01 \\
     				inst\_b20\_3   & 403.057 & 1471.17 & 403.057 & 671.90 & 0.01 \\
     				\hline
     			\end{tabular*}\label{table_cplex_eils}
     		\end{minipage}
     	\end{center}
     \end{table}
     
As shown in Table \ref{table_cplex_eils}, only instances of up to 20 requests are solved using CPLEX. For the rest, none were solved within the limit of 2 hours. One can notice that for these instances, the objective function value of E-ILS is equivalent to that of CPLEX. However, the proposed method requires more computational time.

 \begin{table}[!ht]
     	\begin{center}
     		\begin{minipage}{\textwidth}
     			\caption{Comparison between E-ILS and CPLEX for DARPDP instances based on the three KPIs}	
     			\begin{tabular*}{\textwidth}{@{\extracolsep{\fill}}lcccccc@{\extracolsep{\fill}}}
     				\hline
     				& \multicolumn{3}{@{}c}{CPLEX} & \multicolumn{3}{@{}c@{}}{E-ILS}\\
                         \cline{2-4}\cline{5-7}
     				Instance & {$KPI_{1}[\%]$} & {$KPI_{2}[\%]$} & {$KPI_{3}[min]$} & {$KPI_{1}[\%]$} & {$KPI_{2} [\%]$} & {$KPI_{3}[min]$}\\
     				\hline
     				inst\_a10\_1   & 80  & 50.25   & 6.935   & 80 & 53.69  & 5.992\\
     				inst\_a10\_2   & 80  & 61.87   & 7.960   & 80 & 46.43  & 7.651\\
     				inst\_a15\_1   & 80  & 39.63   & 4.509   & 80 & 29.29  & 4.410\\		
     				inst\_a15\_2   & 80  & 66.88   & 5.594   & 80 & 67.50  & 4.511\\	
     				inst\_a20\_2   & 80  & 36.25   & 7.325   & 80 & 31.25  & 7.295\\
     				inst\_a20\_3   & 80  & 66.12   & 6.416   & 80 & 62.12  & 5.409\\
                        inst\_b10\_1   & 80  & 30.89   & 6.382   & 80 & 13.95  & 2.489\\
     				inst\_b10\_2   & 80  & 50.70   & 11.081  & 80 &  2.10  & 10.181\\	
     				inst\_b15\_1   & 80  & 23.36   & 7.640   & 80 & 31.87  & 5.759\\
     				inst\_b15\_2   & 80  & 25.07   & 7.452   & 80 & 32.99  & 7.293\\	
     				inst\_b20\_2   & 80  & 23.61   & 10.768  & 80 & 31.47  & 10.028\\
     				inst\_b20\_3   & 80  & 47.54   & 9.949   & 80 &  1.81  & 8.093\\		
     				\hline
     			\end{tabular*}\label{Cplex_eils_kpis}
     		\end{minipage}
     	\end{center} 
     \end{table}
 
From table \ref{Cplex_eils_kpis}, it can be seen that the coverage rate is almost equivalent. Actually, the number of served requests is the strict minimum that verifies the $\epsilon-constraint$. For the columns, empty mileage, and ride time, the performance of the E-ILS method is better than CPLEX. Indeed, the empty mileage is about $34\%$ on average for the E-ILS method versus $44\%$ on average for CPLEX. The excess ride time is about 6.59 minutes per client for the E-ILS versus 7.67 minutes per client for the CPLEX.
It should be noted that these KPIs were calculated based on the resulting solutions. Their incorporation within the objective function implies more complexity, requiring dedicated investigation in a further work.

\subsection{Impact of method design components}\label{subsect_9}

In our hybrid approach, several components are applied in order to measure their impact; four combinations are selected and tested as follows:
\begin{itemize}
    \item M0:  the classical ILS uses as a local search phase a simple local search that consists in applying three moves (Relocate, Exchange, R-4-opt) one after the other. The perturbation is based on the same simple local search with the size of the moves increasing by one. And a simple acceptance criterion considering that only the best solutions will be accepted.
    \item M1: Same as above, but use the proposed perturbation.
    \item M2: Same as above, but add a PR phase to intensify the search.
    \item E-ILS: Represents our hybrid approach following the general structure of Algorithm \ref{alg: E-ILS}.
\end{itemize}

In order to compare our method to existing algorithms, another combination of the ILS algorithm with a set partitioning problem is tested which is the ILS-SP method by Malheiros et al.\cite{bib10}. It is important to mention that the algorithm was designed to solve the multi-depot heterogeneous variant of DARP (MDHDARP). An adjusted reimplementation in Python was required to fit our model and data. Results of initial tests show unreasonably long execution times. Thus an additional stopping criterion to the ILS-SP was added; i.e., the same time limits as used in our approach.

\begin{sidewaystable}
     \scriptsize
	\centering
	\caption{Comparison of E-ILS with the three alternate versions on a-type instances}\label{tab6_M}
	\begin{tabular*}{\textwidth}
		{@{\extracolsep{\fill}}lcccccccccccccc@{\extracolsep{\fill}}}
		\hline%
		& \multicolumn{2}{@{}c}{M0} & \multicolumn{3}{@{}c@{}}{M1} & 
		\multicolumn{3}{@{}c@{}}{M2} & \multicolumn{3}{@{}c@{}}{E-ILS} & \multicolumn{3}{@{}c@{}}{ILS-SP} \\
		
		\cline{2-3}\cline{4-6}\cline{7-9}\cline{10-12}\cline{13-15}%
		Instance        & {Best}    & {CPU [s]} &   {Best}  & {CPU [s]} &{gap [\%]}& {Best} & {CPU [s]} &{gap [\%]}& {Best} & {CPU [s]} & {gap [\%]} & {Best} & {CPU [s]} & {gap [\%]}\\
		\hline
		inst\_a10\_1   	&	411.31	&	4.43	&	331.893	& 9.37	    & -19	&	331.893	&	6.97	&  -19    & 331.893	&	1.19    &  -19  & 331.893	&	3.85    &  -19\\
		inst\_a10\_2   	&	222.298	&	3.21	&	151.044	& 19.84	    & -32	&	148.323	&	25.24	&  -33    & 137.1	&	16.15   &  -38  & 137.1	    &	21.59   &  -38\\		
		inst\_a15\_1   	&	409.707	&	41.26   &	392.552	& 56.04	    & -4	&	395.814	&	36.64	&  -3     & 272.922	&	27.97   &  -33  & 421.999	&	35.33   &   3\\		
		inst\_a15\_2  	&	339.154	&	28.31	&	335.874	& 23.04	    & -1	&	318.586	&	45.73	&  -6     & 318.586	&	17.81   &  -6   & 318.586	&	19.25   &  -6\\	
		inst\_a20\_2   	&	567.005	&	126.04	&	509.99	& 124.67	& -10	&	501.013	&	165.10	&  -12    & 410.026	&	93.50   &  -28 & 449.678	&	97.45   &  -21\\
		inst\_a20\_3   	&	383.589	&	79.31	&	377.407	& 309.23	& -2	&	346.548	&	280.42	&  -5     & 353.498	&	136.26  &  -8  & 358.531	&	143.00  &  -7\\
		inst\_a30\_3   	&	996.256 &	84.69	&	989.557	& 272.76	& -1	&	913.185	&	277.15	&  -8     & 909.882	&	446.44  &  -9  & 867.282	&	455.09  &  -13\\
		inst\_a30\_4   	&	983.364 &	548.20	&	947.188	& 731.71	& -4	&	899.68	&	556.93	&  -9     & 812.369	&	612.28  &  -17 & 906.635	&	615.02  &  -8\\	
		inst\_a50\_5  	&	1273.956&	670.96	&	1248.764 & 797.20	& -2	&	1232.311&	444.98	&  -3     & 1147.149&	848.16  & -10 & 1163.03&	853.41  &  -9\\
		inst\_a50\_6   	&	1396.797&	398.45	&	1372.294 & 611.63	& -2	&	1346.459&	486.74	&  -4     & 1271.975&	534.43  &  -9 & 1276.124&	542.46  &  -9\\	
		inst\_a100\_11  &	2725.092&	966.81  &	2451.173 & 1027.44	& -10	&	2329.29	&	786.72	&  -15    & 2209.674&	989.43 &  -19 & 2101.015&	999.37  &  -23\\	
		inst\_a100\_12  &	2367.649&	881.12	&	2190.862 & 1036.31	& -7	&	2045.925&	710.86	& -14  & 1860.426&	961.07  &  -21& 1765.178&	965.10  &  -25\\
		inst\_a200\_22  &	4433.232&	1084.73 &	4234.225 & 1073.47	& -4	&	4227.707&	1102.64	& -5   & 4128.984&	1007.03 &  -7 & 4112.832&	1013.57 &  -7\\
		inst\_a200\_22  &	5261.09 &	1026.43 &	4741.591 & 1105.94	& -10	&	4643.92	&	1124.31	& -12  & 4600.176&	1061.63 &  -13& 4557.271&	1070.96 &  -13\\
		\hline
		Average  & & 424.57 & & 514.19	& & & 432.17 & &  &	482.38 &  &  & 488.25&\\
		CPU[s]   & & & & & & & & & & & & & &\\
		\hline
		Average & &	& &  &  -8	& &	 &	-10   & & &  -17 & & &-14\\
		Gap[\%] & & & & & & & & & & & & & &\\
		\hline
     	\end{tabular*}
     
	\vspace{2\baselineskip}
	\centering
	\caption{Comparison of E-ILS with the three alternate versions on b-type instances}\label{tab7_M}
	
		\begin{tabular*}{\textwidth}
			{@{\extracolsep{\fill}}lcccccccccccccc@{\extracolsep{\fill}}}
			\hline%
			& \multicolumn{2}{@{}c}{M0} & \multicolumn{3}{@{}c@{}}{M1} & 
			\multicolumn{3}{@{}c@{}}{M2} & \multicolumn{3}{@{}c@{}}{E-ILS} & \multicolumn{3}{@{}c@{}}{ILS-SP} \\
			\cline{2-3}\cline{4-6}\cline{7-9}\cline{10-12}\cline{13-15}%
			Instance        & {Best}    & {CPU [s]} &   {Best}  & {CPU [s]} &{gap [\%]}& {Best} & {CPU [s]} &{gap [\%]}& {Best} & {CPU [s]} & {gap [\%]} & {Best} & {CPU [s]} & {gap [\%]}\\
			\hline
			inst\_b10\_1  &	220.093	&	5.38	&	172.758	&	19.57	&	-22	&	172.758	&	17.43	&	-22	&	172.758	&	4.47	&	-22	&	184.904	&	14.30	&	-16	\\
			inst\_b10\_2  &	306.107	&	5.09	&	267.847	&	49.11	&	-12	&	197.29	&	37.34	&	-36	&	197.29	&	28.14	&	-36	&	239.929	&	50.87	&	-22	\\
			inst\_b15\_1  &	345.145	&	24.45	&	328.937	&	27.57	&	-5	&	271.757	&	50.72	&	-21	&	271.757	&	29.50	&	-21	&	286.355	&	49.41	&	-17	\\
			inst\_b15\_2  &	465.27	&	42.61	&	419.624	&	130.54	&	-10	&	398.337	&	137.35	&	-14	&	333.858	&	120.82	&	-28	&	366.005	&	126.76	&	-21	\\
			inst\_b20\_2  &	444.35	&	271.41	&	387.863	&	420.69	&	-13	&	357.538	&	277.99	&	-20	&	253.495	&	349.72	&	-43	&	348.395	&	393.47	&	-22	\\
			inst\_b20\_3  &	639.132	&	455.20	&	582.421	&	265.56	&	-9	&	566.4	&	329.77	&	-11	&	403.057	&	671.90	&	-37	&	471.102	&	678.64	&	-26	\\
			inst\_b30\_3  &	1029.23	&	349.15	&	876.739	&	896.57	&	-15	&	862.052	&	835.43	&	-16	&	784.352	&	871.79	&	-24	&	856.572	&	878.90	&	-17	\\
			inst\_b30\_4  &	844.618	&	724.06	&	780.735	&	625.66	&	-8	&	823.726	&	599.48	&	-2	&	732.154	&	629.29	&	-13	&	816.175	&	690.89	&	-3	\\
			inst\_b50\_5  &	1318.057&	636.38	&	1371.513&	859.22	&	4	&	1296.419&	855.67	&	-2	&	1210.899&	605.95	&	-8	&	1216.253&	664.10	&	-8	\\
			inst\_b50\_6  &	869.341	&	655.65	&	918.229	&	871.46	&	6	&	862.72	&	739.01	&	-1	&	767.605	&	671.15	&	-12	&	807.239	&	702.29	&	-7	\\
			inst\_b100\_11&	2625.673&	919.27	&	2786.48	&	1099.71	&	6	&	2623.454&	930.37	&	0	&	2580.933&	1114.62	&	-2	&	2569.959&	1043.03	&	-2	\\
			inst\_b100\_12&	2292.097&	784.30	&	2143.297&	935.25	&	-6	&	2110.755&	861.50	&	-8	&	1993.452&	1021.65	&	-13	&	1910.7	&	1065.75	&	-17	\\
			inst\_b200\_22&	4591.072&	1041.14	&	4676.657&	1107.45	&	2	&	4696.947&	1041.07	&	2	&	4423.672&	1104.43	&	-4	&	4321.651&	1090.97	&	-6	\\
			inst\_b200\_23&	4073.375&	1131.56	&	4008.19	&	1141.60	&	-2	&	4077.901&	1049.11	&	0	&	3733.878&	1089.74	&	-8	&	3381.286&	1076.08	&	-17	\\
			
			\hline
			Average  & & 503.26 & & 603.57	& & & 554.45 & &  &	593.79 &  &  & 608.96&\\
			CPU[s]   & & & & & & & & & & & & & &\\
			\hline
			Average & &	& &  &  -6	& &	 &	-11   & & &  -19 & & &-14\\
			Gap[\%] & & & & & & & & & & & & & &\\
			\hline
		\end{tabular*}
\end{sidewaystable}

\begin{sidewaystable}
      \scriptsize
	\centering
	\caption{Comparison between E-ILS and the three alternate versions for a-type DARPDP instances based on the three KPIs}\label{tab8}
	\begin{tabular*}{\textheight}{@{\extracolsep{\fill}}lccccccccccccccc@{\extracolsep{\fill}}}
		\hline%
		& \multicolumn{3}{@{}c@{}}{M0} & \multicolumn{3}{@{}c@{}}{M1}& 
		\multicolumn{3}{@{}c@{}}{M2}& \multicolumn{3}{@{}c@{}}{E-ILS} & \multicolumn{3}{@{}c@{}}{ILS-SP}  \\
		\cline{2-4}\cline{5-7}\cline{8-10}\cline{11-13}\cline{14-16}%
		Instance & {$KPI_{1}$} & {$KPI_{2}$} & {$KPI_{3}$} & {$KPI_{1}$} & {$KPI_{2}$} & {$KPI_{3}$} &{$KPI_{1}$} & {$KPI_{2}$} & {$KPI_{3}$} & {$KPI_{1}$} & {$KPI_{2}$} & {$KPI_{3}$} & {$KPI_{1}$} & {$KPI_{2}$} & {$KPI_{3}$}\\
		\hline
		inst\_a10\_1   &  90 & 44 & 5.32  & 80  & 54  & 5.74   & 80 & 54 & 5.92 & 80 & 54 & 5.99 & 80 & 54 & 5.99\\
		inst\_a10\_2   &  80 & 53 & 6.9   & 80  & 62  & 7.78   & 80 & 48 & 7.08 & 80 & 46 & 7.56 & 80 & 46 & 7.56\\
		inst\_a15\_1   &  93 & 32 & 5.17  & 87  & 39  & 5.23   & 87 & 43 & 4.59 & 80 & 29 & 4.41 & 93 & 32 & 5.47 \\
		inst\_a15\_2   &  80 & 67 & 4.73  & 80  & 52  & 4.27   & 80 & 68 & 4.39 & 80 & 68 & 4.51 & 80 & 68 & 4.39\\
		inst\_a20\_2   &  80 & 48 & 7.51  & 80  & 47  & 7.93   & 80 & 43 & 7.40 & 80 & 31 & 7.30 & 80 & 24 & 8.85 \\	
		inst\_a20\_3   &  80 & 56 & 5.30  & 80  & 51  & 6.23   & 80 & 58 & 5.07 & 80 & 62 & 5.41 & 80 & 62 & 5.39 \\		
		inst\_a30\_3   &  80 & 55 & 4.69  & 80  & 49  & 5.17   & 80 & 51 & 5.02 & 80 & 47 & 5.39 & 80 & 48 & 5.19 \\	
		inst\_a30\_4   &  80 & 55 & 6.63  & 80  & 38  & 8.64   & 93 & 30 & 8.85 & 93 & 28 & 8.11 & 83 & 42 & 7.72 \\
		inst\_a50\_5   & 100 & 21 & 8.10  & 98  & 19  & 7.70   & 98 & 20 & 7.78 & 86 & 28 & 8.52 & 90 & 30 & 9.05 \\	
		inst\_a50\_6   &  94 & 21 & 6.61  & 96  & 28  & 7.12   & 94 & 43 & 6.59 & 94 & 44 & 7.25 & 80 & 43 & 7.14 \\			
		inst\_a100\_11 & 80  & 45 & 6.33  & 80  & 43  & 6.49   & 82 & 45 & 5.77 & 80 & 42 & 6.61 & 82 & 42 & 6.37 \\		
		inst\_a100\_12 & 88  & 50 & 6.15  & 80  & 46  & 6.82   & 81 & 39 & 7.24 & 80 & 46 & 7.19 & 80 & 47 & 6.96 \\			
		inst\_a200\_22 & 89  & 43 & 6.75  & 89  & 42  & 7.09   & 89 & 41 & 7.04 & 89 & 46 & 6.48 & 89 & 44 & 6.81 \\	
		inst\_a200\_23 & 89  & 40 & 7.97  & 89  & 44  & 7.03   & 89 & 42 & 7.39 & 89 & 41 & 7.71 & 89 & 42 & 7.26 \\
		\hline
		Average $KPI_{1}[\%]$    & 86 &	  &	  &   84 &   	&	 &	85 &	   &    &	84 &   &  &	83 \\
		\hline
		Average $KPI_{2}[\%]$    &  &	45 &	  &   &  44  	&	 &	 &	 45  &    &	 &  44  &   &	 &  45  &\\
		\hline
		Average $KPI_{3}[min]$   &  &  &	6.30  &   &   &	 6.66 &	  &	   &  6.44  &   &   & 6.60  &   &   & 6.73\\
		\hline
	\end{tabular*} 	

 	\vspace{2\baselineskip}
 	\centering
 	\caption{Comparison between E-ILS and the three alternate versions for b-type DARPDP instances based on the three KPIs}\label{tab9}
 	\begin{tabular*}{\textheight}{@{\extracolsep{\fill}}lccccccccccccccc@{\extracolsep{\fill}}}
 		\hline%
 		& \multicolumn{3}{@{}c@{}}{M0} & \multicolumn{3}{@{}c@{}}{M1}& 
 		\multicolumn{3}{@{}c@{}}{M2} & \multicolumn{3}{@{}c@{}}{E-ILS} & \multicolumn{3}{@{}c@{}}{ILS-SP}\\    
 		\cline{2-4}\cline{5-7}\cline{8-10}\cline{11-13}\cline{14-16}%
 		Instance & {$KPI_{1}$} & {$KPI_{2}$} & {$KPI_{3}$} & {$KPI_{1}$} & {$KPI_{2}$} & {$KPI_{3}$} &{$KPI_{1}$} & {$KPI_{2}$} & {$KPI_{3}$} & {$KPI_{1}$} & {$KPI_{2}$} & {$KPI_{3}$} & {$KPI_{1}$} & {$KPI_{2}$} & {$KPI_{3}$}\\
 		\hline
 		inst\_b10\_1   & 90	  &	61 & 2.40 &	80	& 14 & 2.49	& 80 &	14 & 2.49	& 80 &	14 & 2.49  & 80 &  14 & 5.61\\	
 		inst\_b10\_2   & 80	  &	21 & 8.36 &	80	& 40 & 7.46	& 80 &	2  & 10.18	& 80 &	2  & 10.18 & 80 &  12 & 9.91\\    
 		inst\_b15\_1   & 87	  &	31 & 6.23 &	80	& 31 & 5.35 & 80 &	32 & 5.76	& 80 &	32 & 5.76  & 80 &  32 & 5.71\\    
 		inst\_b15\_2   & 80	  &	42 & 6.65 &	87	& 20 & 7.19 & 87 &	34 & 7.32	& 80 &	33 & 7.29  & 80 &  32 & 7.82\\    
 		inst\_b20\_2   & 95	  &	18 & 9.98 &	80	& 10 & 8.76 & 80 &	33 & 8.59	& 80 &	31 & 10.03 & 80 &  34 & 9.15\\
 		inst\_b20\_3   & 85	  &	14 & 7.46 &	85	& 14 & 9.17	& 90 &	2  & 8.96	& 80 &	2  & 8.09  & 80 &  29 & 6.39\\    
 		inst\_b30\_3   & 90   &	29 & 4.52 &	83	& 55 & 3.92 & 83 &	52 & 4.39	& 83 &	62 & 3.11  & 93 &  51 & 4.12\\
 		inst\_b30\_4   & 90   &	36 & 6.75 &	93	& 50 & 5.74	& 83 &	27 & 7.67	& 93 &	21 & 6.16  & 93 &  33 & 5.79\\   
 		inst\_b50\_5   & 88	  &	34 & 5.94 &	82	& 49 & 5.56 & 84 &	35 & 6.65	& 82 &	37 & 6.16  & 80 &  33 & 6.32\\
 		inst\_b50\_6   & 84	  &	41 & 5.42 &	88	& 57 & 5.00	& 82 &	54 & 5.00	& 80 &	49 & 5.45  & 86 &  41 & 6.18\\    
 		inst\_b100\_11 & 89  &	42 & 6.29 &	89	& 40 & 6.79	& 89   & 46 & 5.88	& 89 &	47 & 6.22 	& 89 &	40 & 6.60\\
 		inst\_b100\_12 & 88  &	47 & 6.12 &	88	& 45 & 5.95	& 88   & 46 & 6.07	& 82 &	44 & 6.65   & 81 &	43 & 5.87\\
 		inst\_b200\_22 & 89  &	48 & 5.88 &	89	& 37 & 6.72	& 89   & 40 & 6.57	& 89 &	40 & 6.63   & 89 &	37 & 7.03\\
 		inst\_b200\_23 & 88.5 &	51 & 5.50 &	80	& 50 & 5.66	& 88.5 & 51 & 5.22	& 80 &	48 & 5.31   & 80 &	48 & 5.11\\
 		\hline
 		Average $KPI_{1}[\%]$       & 87 &	  &	  &   85 &   	&	 &	85 &	   &    &	83 &   &  &	84 \\
 		\hline
 		Average $KPI_{2}[\%]$       &  &	37  &	  &   &  37  	&	 &	 &	 33  &    &	 &  33  &   &	 &  34\\
 		\hline
 		Average $KPI_{3}[min]$      &  &  &	6.25  &   &   &	 6.13 &	  &	   &  6.48  &   &   & 6.40 &   &   & 6.54 \\
 		\hline
 	\end{tabular*}
 \end{sidewaystable}

Tables \ref{tab6_M} and \ref{tab7_M} show the results achieved by the five heuristics for the a-type (b-type) instances, respectively. The percentage gaps between the basic ILS (M0) and the other methods (E-ILS, M1, M2, ILS-SP) are calculated as $100(Best_{Method} - Best_{M0})/Best_{M0}$.

It is worth emphasizing that E-ILS and ILS-SP significantly outperforms the basic ILS (M0) with an average of 17\%(19\%) and 14\%(14\%) respectively increase in solution quality for all a-type (b-type) DARPDP instances, respectively. However, this comes at the expense of more computational time but is still reasonable, as indicated in tables \ref{tab6_M} and \ref{tab7_M}. Besides, the solutions obtained by $M1$ with LNS-based perturbation are 8\% (6\%) on average better than $M0$, which is already, an important improvement, and hence this approach effectively diversify the search space. As one might expect, adding a PR phase intensify deeply the search and gives an improvement of 10\% (11\%) compared to M0. It is worth mentioning that the considerable improvement achieved by the E-ILS is related to exchanging the Simple Local Search for a Learning Local Search.
One can notice that on large instances, the ILS-SP performs slightly better. However, for the small instances our proposed approach is much more effective.

In the comparative tables \ref{tab8} and \ref{tab9}, one can notice that there is no assertion that one method is the most efficient on the three key performance indicators simultaneously. Although we can draw the following conclusions.

\begin{figure}[h!]
     \begin{center}
         \includegraphics[width=0.5\textwidth,height=0.3\textwidth]{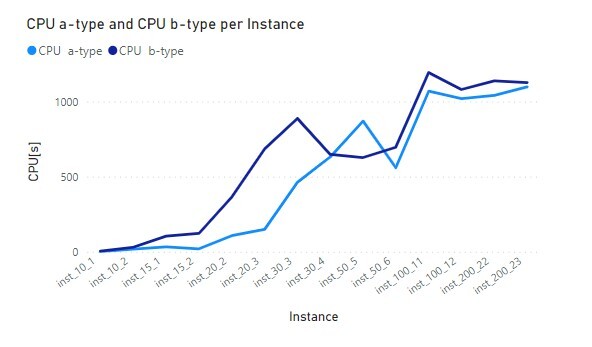}
         \caption{\textbf{Running time per instance type} The figure illustrates The variation of execution time depending on the types of instances}\label{fig2}
    \end{center}
\end{figure}

According to figure \ref{fig2}, it can be seen that the computation time obviously increases with increasing the number of requests. Comparing the scenarios, a-type, and b-type instances. As one might expect, the complexity of the problem increases with b-type instances, which is consistent with the reported increase in resolution time.

\begin{figure}[h!]
     \begin{center}
         \includegraphics[width=0.96\textwidth,height=0.3\textwidth]{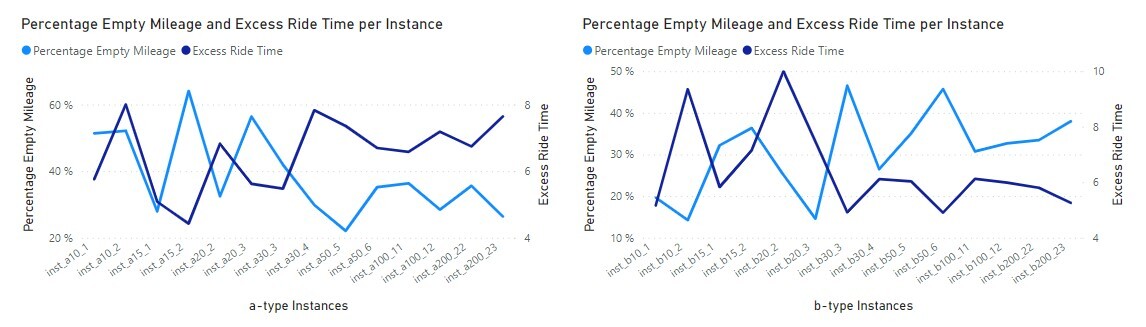}
         \caption{\textbf{Correlation between key performance indicators} The two figures illustrate the average percentage of Empty Mileage and the average Excess Ride Time on the solutions for a-type test instances and b-type test instances}\label{fig3}
         \end{center}	
\end{figure}

Through the analysis of the figure \ref{fig3}, one can notice that the percentage of empty mileage decreases with the increase in the excess ride time. Therefore, as one might expect, an extra ride time is related to an indirect trip. Interestingly, this decrease in direct trips also comes with a general reduction of the scenario with one request onboard the vehicle.

\section{Conclusion}\label{Sec6}
\quad In this research, the dial-a-ride problem was adapted to handle a more practical application. The new problem considers driver preferences in terms of working during times and in areas that suit their personal preferences. To ensure this, we have added for each driver individual origin, destination and time window to reach his destination. The driver’s preferences modeled by additional threshold constraint which define the maximum allowable tolerance time for reaching the driver’s destination, and this deviation is minimised further in the objective function. We proposed a new efficient E-ILS algorithm where the advantages of multiple heuristics are fully investigated and combined into a unified ILS framework, including learning local search, path-relinking, and large neighborhood search. The computational results, using generated instances based on real-world data and analyzing three key indicator performances, confirm that the metaheuristic approach is more efficient than CPLEX at producing high-quality solutions in a reasonable time for the small DARPDP instances and outperforms the basic ILS for the large ones.\\
\quad Our future work will focus on three aspects: (1) Addressing the online version of the proposed problem is an interesting research direction. (2) The problem is challenging to resolve because of the multiplicity and complexity of the constraints. The performance of the proposed model can be improved either by adding cuts to relax hard constraints or by decomposing the problem into sub-problems. (3) We paid considerable attention to designing the algorithm; it seems to be a very promising approach. However, the results of preliminary experiments show that we still need more tests on the parameters besides the calibration of probabilities to call different method components. Hence, additional verification could result in better solution performance and convergence.

\bibliographystyle{unsrtnat}
\bibliography{references}  

\end{document}